\documentclass[reqno,draft,a4paper]{amsart}
\usepackage{times}
\usepackage{mathrsfs}
\usepackage{latexsym,amsopn,amssymb,amsmath}
\usepackage{amsthm}
\usepackage{amsfonts,amsbsy,amscd,stmaryrd}
\usepackage{paralist}
\usepackage{geometry}
\oddsidemargin %-0.1cm
\evensidemargin %-0.1cm
\newcommand{\mbc}{\mathbb{C}}

\newcommand{\mbr}{\mathbb{R}}

\newcommand{\cc}{\mathcal{C}}

\newcommand{\cf}{\mathcal{F}}
\newcommand{\cg}{\mathcal{G}}
\newcommand{\ch}{\mathcal{H}}

\def\rmd{\text{d}}
\def\rme{\mathrm{e}}
\def\rmo{\mathrm{o}}

\newcommand{\Co}{\cc_{\rmo}}

\renewcommand{\proof}{\noindent{\bf Proof: }}

\def\braket#1#2{\langle{#1}|{#2}\rangle}

\def\jap#1{\langle {#1} \rangle}

\def\slim{\mbox{\rm s-}\!\lim}
\def\wlim{\mbox{\rm w-}\!\lim}

\def\qed{\hfill \raisebox{0.5ex}{\framebox[1.6ex]{
                                       \rule[0ex]{0ex}{0.3ex} }}}

\def\build#1_#2^#3{\mathrel{\mathop{\kern 0pt#1}\limits_{#2}^{#3}}}
\newcounter{PAR}[section]
\def\thePAR{\arabic{PAR}.}
\def\PAR{\par\addvspace{2 ex}
\noindent 
\refstepcounter{PAR}{\bf \thePAR}
}
\newtheorem{theorem}{Theorem}%[section]
\newtheorem{lemma}[theorem]{Lemma}
\newtheorem{proposition}[theorem]{Proposition}

\newtheorem{remark}[theorem]{\bf Remark}
\newtheorem*{acknowledgment}{Acknowledgment}
\long\def\symbolfootnote[#1]#2{\begingroup%
\def\thefootnote{\fnsymbol{footnote}}\footnote[#1]{#2}\endgroup} 
\begin{document}

\title{Hamiltonians with purely discrete spectrum} 

\author[Vladimir GEORGESCU]{Vladimir GEORGESCU} \address{CNRS and
  University of Cergy-Pontoise 95000 Cergy-Pontoise, France}
\email{vlad@math.cnrs.fr} 

\date{\today} 

\thanks{}

\begin{abstract}
We discuss criteria for a self-adjoint operator on $L^2(X)$ to have
empty essential spectrum. We state a general result for the case of a
locally compact abelian group $X$ and give examples for $X=\mbr^n$.
\end{abstract}

\maketitle 

\vspace{-5mm}

\PAR
Let $\Delta$ be the positive Laplacian on $\mbr^n$.  We set
$B_a(r)=\{x\in\mbr^n \mid |x-a|\leq r\}$ and $B_a=B_a(1)$. 

\begin{proposition}\label{pr:lap}
Let $V$ be a real locally integrable function on $\mbr^n$ such that:

\begin{compactenum}
\item[{\rm(i)}] if $\lambda >0$ then the measure $\omega_\lambda(a)$
  of the set $\{x\in B_a \mid V(x)<\lambda \}$ satisfies
  $\lim_{a\to\infty}\omega_\lambda(a)=0$,
\item[{\rm(ii)}] the negative part of $V$ satisfies $V_-\leq
  \mu\Delta+\nu$ for some positive real numbers $\mu,\nu$ with
  $0<\mu<1$.
\end{compactenum}
Then the spectrum of the self-adjoint operator $H$ associated to the
form sum $\Delta+V$ is purely discrete.
\end{proposition} 

\begin{remark}\label{re:lap1}{\rm
Let $V_\pm=\max\{\pm V, 0\}$ and for each $\lambda>0$ let
$\Omega_\lambda=\{ x \mid V_+(x)<\lambda \}$. Then
$\omega_\lambda(a)$ is the measure of the set
$B_a\cap\Omega_\lambda$.  From Lemma \ref{lm:wv} it follows that the
condition (i) is equivalent to
\begin{equation}\label{eq:wdecay}
\lim_{a\to\infty}\int_{B_a}\frac{\rmd x}{1+V_+(x)} =0.
\end{equation}
}\end{remark}

\begin{remark}\label{re:lap2}{\rm
From Lemma \ref{lm:omega} we get
$\lim_{a\to\infty}\omega_\lambda(a)=0$ if
$\int_{\Omega_\lambda}\omega_\lambda^p \rmd x < \infty $ for some
$p>0$.  Thus Theorems 1 and 3 from \cite{S} are consequences of
Proposition \ref{pr:lap}.  In the case $V\geq0$ Proposition
\ref{pr:lap} is a consequence of Theorem 2.2 from \cite{MS}.  More
general results will be obtained below.  Note, however, that our
techniques are not applicable in the framework considered in Theorem
2 from \cite{S} and in \cite{WW}.  }\end{remark}

Proposition \ref{pr:lap} is very easy to prove if condition
\eqref{eq:wdecay} is replaced by $\lim_{x\to\infty}V_+(x)=\infty$.
In fact, let us consider an arbitrary locally compact space $X$ and
let $\ch$ be a Hilbert $X$-module, i.e. $\ch$ is a Hilbert space and
a nondegenerate $*$-morphism $\phi\mapsto\phi(Q)$ of $\Co(X)$ into
$B(\ch)$ is given. For example, one may take $\ch=L^2(X,\mu)$ for
some Radon measure $\mu$. Then we have the following simple
compactness criterion: \emph{if $R$ is a bounded self-adjoint
  operator on $\ch$ such that (i) if $\phi\in\Co(X)$ then $\phi(Q)R$
  is a compact operator, (ii) one has $\pm R\leq\theta(Q)$ for some
  $\theta\in\Co(X)$, then $R$ is a compact operator}. Indeed, note
first that the operator $R\phi\equiv R\phi(Q)$ will also be compact
for all $\phi\in\Co(X)$. Let $\varepsilon>0$ and let us choose
$\phi$ such that $0\leq\phi\leq1$ and
$\theta\phi^\perp\leq\varepsilon$, where $\phi^\perp=1-\phi$. Then
$\pm\phi^\perp R\phi^\perp \leq \phi^\perp \theta
\phi^\perp\leq\varepsilon$ which implies $\|\phi^\perp
R\phi^\perp\|\leq\varepsilon$. So we have $\|R-\phi R-\phi^\perp
R\phi\|\leq\varepsilon$ and $\phi R+\phi^\perp R\phi$ is a compact
operator.  Now let us say that a self-adjoint operator $H$ on $\ch$
is \emph{locally compact} if $\phi(Q)(H+i)^{-1}$ is compact for all
$\phi\in\Co(X)$. Then we get: \emph{ If $H$ is a locally compact
  self-adjoint operator on $\ch$ and if there is a continuous
  function $\Theta:X\to\mbr$ such that $\lim_{x\to\infty}
  \Theta(x)=+\infty$ and $H\geq\Theta(Q)$, then the spectrum of $H$
  is purely discrete} (the nondegeneracy of the morphism is needed
for the definition of $\Theta(Q)$ for unbounded $\Theta$).

\PAR
On the other hand, Proposition \ref{pr:lap} can be significantly
generalized. For example, $\Delta$ may be replaced by a higher order
operator with matrix valued coefficients and $V$ does not have to be
a function. These results are consequences of the following
``abstract'' fact. We fix a locally compact abelian group $X$,
choose a finite dimensional Hilbert space $E$, and define
$\ch=L^2(X)\otimes E$.  For $a\in X$ and $k\in X^*$ (the dual
locally compact abelian group) we denote $U_a$ and $V_k$ the unitary
operators on $\ch$ given by
\[
(U_af)(x)=f(x+a) \hspace{2mm} \text{and}\hspace{2mm} 
(V_kf)(x)=k(x)f(x).
\]
We denote additively the operations both in $X$ and in $X^*$ and
denote $0$ their neutral elements.

\begin{theorem}\label{th:main}
Let $H$ be a self-adjoint operator on $\ch$ such that for some
(hence for all) $z\in\mbc$ not in the spectrum of $H$ the
operator $R=(H-z)^{-1}$ satisfies
\begin{equation}\label{eq:main}
\lim_{k\to0} \|V_kR V_k^*-R\|=0, \quad \lim_{a\to0} \|(U_a-1)R\|=0. 
\end{equation}
Then $H$ has purely discrete spectrum if and only if
$\wlim_{a\to\infty} U_aRU_a^*=0$. 
\end{theorem}

\proof If the spectrum of $H$ is purely discrete then $R$ is compact
so $\wlim_{a\to\infty} U_aRU_a^*=0$. The reciprocal assertion follows
from Theorem 1.2 from \cite{GI}. Indeed, with the terminology used
there, all the localizations at infinity of $H$ will be equal to
$\infty$ hence the essential spectrum of $H$ will be empty.  \qed

Some notations: if $\phi$ is a $B(E)$-valued Borel function on $X$
then $\phi(Q)$ is the operator of multiplication by $\phi$ on $\ch$;
if $\psi$ is a similar function on $X^*$ then $\psi(P)=\cf^{-1} M_\psi
\cf$, where $\cf$ is the Fourier transformation and $M_\psi$ is the
operator of multiplication by $\psi$ on $L^2(X^*)\otimes E$. Note that
$V_k\psi(P)V_k^*=\psi(P+k)$.

If $\phi\in L^\infty(X)$ and $\phi\geq0$ then it is easy to check
that $\wlim_{a\to\infty}U_a\phi(Q)U_a^*=0$ if and only if
$\slim_{a\to\infty}\phi(Q)U_a=0$ and also if and only if there is a
compact neighborhood of the origin $W$ such that
$\lim_{a\to\infty}\int_{a+W}\phi\rmd x=0$. Then we say that $\phi$
is \emph{weakly vanishing (at infinity)}. See Section 6 in \cite{GG}
for further properties of this class of functions. Below $W$ is a
compact neighborhood of the origin, $W_a=a+W$, and we denote $|M|$
the Haar measure of a set $M$.

\begin{lemma}\label{lm:wv}
A positive function $\phi\in L^\infty(X)$ is weakly vanishing if and
only if for any number $\lambda>0$ the set $\Omega^\lambda=\{ x \mid
\phi(x)>\lambda \}$ has the property
$\lim_{a\to\infty}|W_a\cap\Omega^\lambda|=0$.
\end{lemma}
This follows from the estimates
\[
\lambda|W_a\cap\Omega^\lambda|\leq \int_{W_a}\phi\, \rmd x\leq 
\|\phi\|_{L^\infty} |W_a\cap\Omega^\lambda| + \lambda |W|.
\]

\begin{proposition}\label{pr:main}
Let $H$ be an invertible self-adjoint operator satisfying
\eqref{eq:main} and such that $\pm H^{-1}\leq \phi(Q)$ for some
weakly vanishing function $\phi$. Then $H$ has purely discrete
spectrum.
\end{proposition}
Indeed, we may take $R=H^{-1}$ and then for any $f\in\ch$ we have
$|\braket{f}{U_aRU_a^*}|\leq |\braket{f}{U_a\phi(Q)U_a^*}|$.

{\bf Proof of Proposition \ref{pr:lap}:} Here $X=\mbr^n$ and we
identify as usual $X$ with its dual by setting $k(x)=\rme^{i kx}$
for $x,k\in X$. Then if $P_j=-i\partial_j$ and $P=(P_1,\dots,P_n)$
we get $V_kPV_k^*=P+k$. To simplify notations we write $H$ for
$\Delta+V+1+\nu$, so that $H\geq(1-\mu)\Delta+V_++1\geq
V_++1\geq1$. Then observe that the form domain of $H$ is $\cg\equiv
D(H^{1/2})=\{f\in\ch^1 \mid V_+^{1/2}f\in L^2\}$ where $\ch^1$ is
the first order Sobolev space. Thus $R=H^{-1}:L^2\to\ch^1$ is
continuous and this implies the second part of condition
\eqref{eq:main}.  On the other hand, $H$ extends to a continuous
bijective operator $\cg\to\cg^*$ whose inverse is an extension of
$R$ to a continuous map $\cg^*\to\cg$. We keep the notations $H,R$
for these extensions.  Clearly $V_k$ leaves invariant $\cg$ hence
extends to a continuous operator on $\cg^*$ and the groups of
operators $\{V_k\}$ are of class $C_0$ in both spaces. Now
$H_k:=V_kHV_k^*=(P+k)^2+V=H+2kP+k^2$ in $B(\cg,\cg^*)$
so if $R_k:=V_kRV_k^*$ then
\[
R_k-R=R_k(H-H_k)R= -R_k(2kP+k^2)R
\]
in $B(\cg^*,\cg)$.  Now clearly the first part of \eqref{eq:main} is
fulfilled. Finally, it suffices to show that $H^{-1}\leq\phi(Q)$ for
a weakly vanishing function $\phi$. But $H\geq 1+V_+$ and we may
take $\phi=(1+V_+)^{-1}$ due to \eqref{eq:wdecay}.  \qed

Remark \ref{re:lap2} is a consequence of the next result.

\begin{lemma}\label{lm:omega}
Let $\Omega\subset\mbr^n$ be a Borel set and let
$\omega:\mbr^n\to\mbr$ be defined by $\omega(a)=|B_a\cap\Omega|$.  If
$\omega^p$ is integrable on $\Omega$ for some $p>0$ then
$\omega(a)\to0$ as $a\to\infty$.
\end{lemma}
\proof The main point is the following observation due to Hans Henrik
Rugh: \emph{let $\nu$ be the minimal number of (closed) balls of
  radius $1/2$ needed to cover a ball of radius one; then for any $a$
  there is a Borel set $A_a\subset B_a\subset\Omega$ with
  $|A_a|\geq\omega(a)/\nu$ such that $\omega(x)\geq\omega(a)/\nu$ if
  $x\in A_a$}. Indeed, let $N$ be a set of $\nu$ points such that
$B_a\subset\cup_{b\in N} B_b(1/2)$. If $D_b=B_a\cap B_b(1/2)$ then
$\omega(a)\leq\sum_b|D_b\cap\Omega|$ hence there is $b(a)$ such that
$A_a=D_{b(a)}\cap\Omega$ satisfies $|A_a|\geq\omega(a)/\nu$. Since
$A_a$ has diameter smaller than one, for $x\in A_a$ we have
$A_a\subset B_x\cap\Omega$ hence $\omega(x)\geq|A_a|$, which proves
the remark. Now let us set $R=|a|-1$ and denote $\Omega(R)$ the set of
points $x\in\Omega$ such that $|x|\geq R$. Then we have
\[
\int_{\Omega(R)} \omega^p \rmd x \geq
\int_{A_a} \omega^p \rmd x\geq [\omega(a)/\nu]^{p+1}
\]
which clearly implies the assertion of the lemma.
\qed

\PAR We present here some consequences of Proposition \ref{pr:main}.
We refer to \cite{GI} for general classes of operators verifying
condition \eqref{eq:main} and consider here only some particular
cases. We mention that \emph{if $H$ is a bounded from below operator
  satisfying \eqref{eq:main} and if $\theta:\mbr\to\mbr$ is a
  continuous function such that $\theta(\lambda)\to+\infty$ when
  $\lambda\to+\infty$ the $\theta(H)$ also satisfies \eqref{eq:main}}.

If $R\in B(\ch)$ satisfies the first part of \eqref{eq:main} we say
that $R$ is a \emph{regular operator} (or $Q$-regular). The regularity
of the resolvent of a differential operators on $\mbr^n$ is easy to
check because $V_kPV_k^*=P+k$, cf. the proof of Proposition
\ref{pr:lap}.  The second part of \eqref{eq:main} is equivalent to the
existence of a factorization $R=\psi(P)S$ with $\psi\in\Co(X^*)$ and
$S\in B(\ch)$. If $X=\mbr^n$ then it suffices that the domain of $H$
be included in some Sobolev space $\ch^m$ with $m>0$ real.  We now
give an extension of Proposition \ref{pr:lap} which is proved in
essentially the same way. We assume $X=\mbr^n$ and work with Sobolev
spaces but a similar statement holds for an arbitrary $X$: it suffices
to replace the function $\jap{k}^m$ which defines $\ch^m$ by an
arbitrary weight \cite{GI} and the ball $B_a$ by $a+W$ where $W$ is a
compact neighborhood of the origin.

\begin{proposition}\label{pr:gen}
Let $H_0$ be a bounded from below self-adjoint operator on $\ch$ with
form domain equal to $\ch^m$ for some real $m>0$ and satisfying
$\lim_{k\to0} V_k H_0 V_k^*=H_0$ in norm in $B(\ch^m,\ch^{-m})$.  Let
$V$ be a positive locally integrable function such that
$\lim_{a\to\infty}|\{x\in B_a \mid V(x)<\lambda \}|=0$ for each
$\lambda >0$. Then the self-adjoint operator $H$ associated to the
form sum $H_0+V$ has purely discrete spectrum.
\end{proposition}

Let $h:X\to B(E)$ be a continuous symmetric operator valued function
with $c'|p|^{2m}\leq h(p) \leq c''|p|^{2m}$ (as operators on $E$)
for some constants $c',c''>0$ and all large $p$. Let
$W:\ch^m\to\ch^{-m}$ be a symmetric operator such that 
$W\geq-\mu h(P)-\nu$ with $\mu<1$ and such that $V_k W V_k^*\to
W$ in norm in $B(\ch^m,\ch^{-m})$ as $k\to0$. Then the form sum
$h(P)+W$ is bounded from below and closed on $\ch^m$ and the
self-adjoint operator $H_0$ associated to it satisfies the
conditions of Proposition \ref{pr:gen}.  

Assume that $m\geq1$ is an integer and let
$L=\sum_{\alpha,\beta}P^\alpha
a_{\alpha\beta}(Q)P^\beta:\ch^m\to\ch^{-m}$ where $\alpha,\beta$ are
multi-indices of length $\leq m$ and $a_{\alpha\beta}$ are functions
$X\to B(E)$ such that $a_{\alpha\beta}(Q)$ is a continuous map
$\ch^{m-|\beta|}\to\ch^{|\alpha|-m}$. If $\braket{f}{Lf}\geq \mu
\|f\|^2_{\ch_m} -\nu\|f\|^2_{\ch}$ for some $\mu,\nu>0$ then $L$ is a
closed bounded from below form on $\ch^m$ and the self-adjoint
operator $H_0$ associated to it verifies Proposition \ref{pr:gen}.

{\bf Note added July 2014:} The theory can be extended to metric
spaces by using the $C^*$-algebra introduced and studied in my paper
``On the structure of the essential spectrum of elliptic operators on
metric spaces'', J. Funct. Analysis 260, 1734--1765 (2011) and
arXiv:1003.3454.

\begin{acknowledgment}{\rm
We thank Hans Henrik Rugh for the remark which made Lemma
\ref{lm:omega} obvious.  }\end{acknowledgment}

\end{document}